\documentclass[a4paper,reqno,12pt]{amsart}
\usepackage{verbatim}
\usepackage{amssymb,amsmath,amsthm}
\usepackage{amsfonts}
\usepackage{array}

\newtheorem{theorem}{Theorem}
\newtheorem{lemma}{Lemma}

\theoremstyle{remark}

\newcommand*\FF[4]{{}_{2}F_1\bigg[\begin{matrix}
#1,&#2\\&#3\end{matrix}\,;\,#4\bigg]}

\begin{document}

\setcounter{page}{1}

\title[]{Three pairs of congruences concerning \\sums of central binomial coefficients}

\author{Guo-Shuai Mao}

\address{Department of Mathematics, Nanjing
University of Information Science and Technology, Nanjing 210044,  People's Republic of China}
\email{maogsmath@163.com}

\author{Roberto~Tauraso}

\address{Dipartimento di Matematica, 
Universit\`a di Roma ``Tor Vergata'', 
via della Ricerca Scientifica,
00133 Roma, Italy}
\email{tauraso@mat.uniroma2.it}



\begin{abstract} 
Recently the first author proved a congruence proposed in 2006 by Adamchuk: 
$\sum_{k=1}^{\lfloor\frac{2p}{3}\rfloor}\binom{2k}{k}\equiv 0\pmod{p^2}$ for any prime $p=1 \pmod{3}$.
In this paper, we provide more examples (with proofs) of congruences of the same kind
$$\sum_{k=1}^{\lfloor\frac{ap}{r}\rfloor}\binom{2k}{k}x^k \pmod{p^2}$$
where $p$ is a prime such that $p\equiv 1  \pmod{r}$, $a/r$ is a fraction in $(1/2,1)$ and $x$ is a $p$-adic integer.
The key ingredients are the $p$-adic Gamma functions $\Gamma_p$ and a special class of computer-discovered hypergeometric identities. 
\end{abstract}

\maketitle

\section{Introduction}
In the past decades, many people studied congruences modulo a power of a prime $p$ for sums of binomial coefficients. A certain number  concern  the central binomial coefficients and have the form
$$\sum_{k=0}^n\binom{2k}{k}x^k$$
where the upper limit $n$ is usually $p-1$ or $\frac{p-1}{2}$. We would like to mention a few of them that we will need later. By \cite[p. 156]{MattareiTauraso2013} and \cite[Theorem 1.1]{SunZW2013}, for any prime $p>3$, 
\begin{align}
\label{ncong1a}
&\sum_{k=0}^{p-1}\binom{2k}{k}(-2)^k
\equiv 1-\frac{4pq_p(2)}{3}\pmod{p^3}\\
\label{ncong1b}
&\sum_{k=0}^{\frac{p-1}{2}}\frac{\binom{2k}k}{(-32)^k}
\equiv\left(\frac{2}{p}\right)\left(1+\frac{pq_p(2)}{6}-\frac{p^2q^2_p(2)}{8}\right)\pmod{p^3}
\end{align}
where $q_p(a)=\frac{a^{p-1}-1}{p}$ is the so-called \emph{Fermat quotient}.
Moreover, by \cite[Theorem 1.3]{PanSun2014}, \cite[Theorem 1.1]{SunZW2013}, for any prime $p>5$, 
\begin{align}\label{ncong2a}
&\sum_{k=0}^{p-1}\binom{2k}{k}(-1)^k\equiv
\left(\frac{p}{5}\right)
\left(1-2pf_p\right)\pmod{p^3}\\
\label{ncong2b}
&\sum_{k=0}^{\frac{p-1}{2}}\frac{\binom{2k}{k}}{(-16)^k}\equiv\left(\frac{p}{5}\right)
\left(1+\frac{pf_p}{2}\right)\pmod{p^3}
\end{align}
where $f_p=\frac{F_{p-\left(\frac{p}{5}\right)}}{p}$ is the \emph{Fibonacci quotient} and $F_n$ denotes the $n$th Fibonacci number.

Less common are the examples where the upper limit of the sum is strictly between $\frac{p-1}{2}$ and $p-1$.
In 2014, Pan and Sun \cite{PanSun2014} proved that for any prime $p\equiv1\pmod 4$,
$$\sum_{k=0}^{\lfloor\frac{3p}{4}\rfloor}\frac{\binom{2k}{k}}{(-4)^k}\equiv\left(\frac{2}{p}\right)\pmod{p^2}.$$
Recently Mao \cite{Mao2020a,Mao2020b} proved
that for any prime $p\equiv1\pmod3$, we have
$$\sum_{k=1}^{\lfloor\frac{2p}3\rfloor}\binom{2k}{k}\equiv 0
\quad\text{and}\quad
\sum_{k=0}^{\lfloor\frac{5p}6\rfloor}\frac{\binom{2k}{k}}{16^k}\equiv\left(\frac{3}{p}\right)\pmod{p^2}.
$$
The first congruence has been conjectured in 2006 by Adamchuk \cite{Adamchuk2006}.

The main purpose of this paper is to show more congruences of the same kind. 

\begin{theorem}\label{Thadam}
Let $p$ be a prime. If $p\equiv1\pmod 3$ then
\begin{equation}\label{pair1}
\sum_{k=1}^{\lfloor\frac{2p}3\rfloor}\binom{2k}{k}(-2)^k\equiv0,\qquad
\sum_{k=0}^{\lfloor\frac{5p}6\rfloor}\frac{\binom{2k}k}{(-32)^k}\equiv\left(\frac{2}{p}\right)\pmod{p^2}.
\end{equation}
If $p\equiv1\pmod 5$ then
\begin{equation}\label{pair2}
\sum_{k=1}^{\lfloor\frac{4p}5\rfloor}\binom{2k}{k}(-1)^k\equiv 0,\qquad
\sum_{k=1}^{\lfloor\frac{7p}{10}\rfloor}\frac{\binom{2k}{k}}{(-16)^k}\equiv 0
\pmod{p^2},
\end{equation}
and
\begin{equation}\label{pair3}
\sum_{k=1}^{\lfloor\frac{3p}5\rfloor}\binom{2k}{k}(-1)^k\equiv 0,\qquad
\sum_{k=1}^{\lfloor\frac{9p}{10}\rfloor}\frac{\binom{2k}{k}}{(-16)^k}\equiv 0\pmod{p^2}.
\end{equation}
\end{theorem}
The first pair of congruences seems to have (so far) slipped the attention of the mathematical community whereas the other pairs appeared as conjectures in \cite{SunZW2013}.

The paper is organized as follows. In Section 2 we recall some properties of the $p$-adic Gamma function and some known congruences involving harmonic sums. Section 3 is devoted to a auxiliary binomial sum defined as
$$S_m(x)=\frac{1}{\binom{2m}{m}}\sum_{k=0}^{m-1}
\frac{\binom{2k}{k}x^{k-m}}{k-m}.$$
It turns out that $S_m(x)$ can be evaluated modulo a prime for some special values of $m$ and $x$ by using a class of finite hypergeometric identites computer-discovered by Shalosh B. Ekhad (well-known co-author of Doron Zeilberger).   
In Section 4, we explain how $S_m(x)$ is related to our main results. Finally, we give a proof of Theorem \ref{Thadam} which reveals why the congruences \eqref{pair1}, \eqref{pair2}, and \eqref{pair3} come in pairs.

\section{Notations and preliminaries}

The generalized hypergeometric function is defined as
\begin{equation*}
{}_{m+1}F_m\bigg[\begin{matrix}
\alpha_0&\alpha_1&\ldots&\alpha_m\\
&\beta_1&\ldots&\beta_m
\end{matrix}\bigg|\,z\bigg]=\sum_{k=0}^{\infty}\frac{(\alpha_0)_k(\alpha_1)_k\cdots(\alpha_m)_k}{(\beta_1)_k\cdots(\beta_m)_k}\cdot\frac{z^k}{k!},
\end{equation*}
where $\alpha_0,\ldots,\alpha_m,\beta_1,\ldots,\beta_m,z\in\mathbb{C}$ and
$$
(\alpha)_k=\begin{cases}\alpha(\alpha+1)\cdots(\alpha+k-1) &\text{if }k\geq 1,\\
1 &\text{if }k=0.\end{cases}
$$
is the rising factorial or the Pochhammer
symbol.

The $p$-adic gamma function $\Gamma_p$ is defined as the continuous extension to the ring of all $p$-adic integers $\mathbb{Z}_p$ of the sequence
$$
n \to (-1)^n\prod_{\substack{1\leq j<n\\ (k,p)=1}}k
$$
where $p$ is a prime and $n$ is a positive integer.

It follows that $\Gamma_p(0)=1$ and, for $x\in \mathbb{Z}_p$,
$$
\frac{\Gamma_p(x+1)}{\Gamma_p(x)}=
\begin{cases}-x&\text{if }|x|_p=1,\\
-1&\text{if }|x|_p>1.
\end{cases}
$$
where $|\cdot|_p$ is the $p$-adic norm.
Throughout the whole paper, we will only use the most basic properties of $\Gamma_p$ (see \cite[Chap. 7]{Robert2000} as a general reference). In particular, we recall that (see for example \cite[(3.2)]{		MaoPan2017}) for any $\alpha\in\mathbb{Z}_p$,
\begin{equation}\label{de}
\frac{\Gamma_p'(x)}{\Gamma_p(x)}\equiv \Gamma_p'(0)+H_{p-\langle-x\rangle_p-1}\pmod{p},
\end{equation}
where $H_n=\sum_{k=1}^n\frac1k$ is the $n$th harmonic number and $\langle x\rangle_p$ is the least nonnegative residue of $x$ modulo $p$.

\begin{lemma} \label{LHf}
For any prime $p>3$,
\begin{align*}
&H_{\lfloor \frac{p}{2}\rfloor}\equiv-2q_p(2),\quad 
H_{\lfloor \frac{p}{4}\rfloor}\equiv-3q_p(2)
\pmod{p},\\ 
&H_{\lfloor \frac{p}{3}\rfloor}
\equiv-\frac{3q_p(3)}{2},\quad 
H_{\lfloor \frac{p}{6}\rfloor}
\equiv-2q_p(2)-\frac{3q_p(3)}{2}
\pmod{p}.
\end{align*}
For any prime $p>5$,
\begin{align*}
&H_{\lfloor\frac{p}{5}\rfloor}\equiv 
 -\frac{5q_p(5)}{4}-\frac{5f_p}{4},\quad 
H_{\lfloor\frac{2p}{5}\rfloor}\equiv 
 -\frac{5q_p(5)}{4}+\frac{5f_p}{4}\pmod{p},\\
&H_{\lfloor\frac{p}{10}\rfloor}\equiv -2q_p(2)-\frac{5q_p(5)}{4}
-\frac{15f_p}{4},\quad 
H_{\lfloor\frac{3p}{10}\rfloor}\equiv -2q_p(2)-\frac{5q_p(5)}{4}
+\frac{15f_p}{4}\pmod{p}.
\end{align*}
\end{lemma}
\begin{proof} For a proof of the first group of congruences see \cite[(41)-(44)]{Lehmer1938}. As regards the second group, the first two congruences follow from 
$$H_{\lfloor\frac{2p}{5}\rfloor}+H_{\lfloor\frac{p}{5}\rfloor}\equiv -\frac{5q_p(5)}{2} ,\quad
H_{\lfloor\frac{2p}{5}\rfloor}-H_{\lfloor\frac{p}{5}\rfloor}\equiv \frac{5f_p}{2} \pmod{p}$$
see \cite[(15)]{Lerch1905} and \cite[(4.7)]{Williams1991} respectively. Finally, the last two congruences are obtained by considering two more linear relations
\cite[(12) with $n=10$]{Vandiver2017},
$$10(q_p(2)+q_p(5)) \equiv 10q_p(10)\equiv -\sum_{k=1}^{9} H_{\lfloor\frac{kp}{10}\rfloor}
\equiv -2\sum_{k=1}^{4} H_{\lfloor\frac{kp}{10}\rfloor}-H_{\lfloor\frac{p}{2}\rfloor}\pmod{p},$$ 
and
\begin{align*}
5f_p\equiv 2\sum_{k=1}^{\lfloor\frac{4p}{5}\rfloor}\frac{(-1)^k}{k}
= 2\sum_{k=1}^{\lfloor\frac{2p}{5}\rfloor}\frac{1}{2k}
+2\sum_{k=\lfloor\frac{p}{10}\rfloor+1}^{\frac{p-1}{2}}\frac{1}{2k-p}
\equiv H_{\lfloor\frac{2p}{5}\rfloor}+H_{\lfloor\frac{p}{2}\rfloor}-H_{\lfloor\frac{p}{10}\rfloor}\pmod{p},
\end{align*}
where the congruence on the left is given by \cite[Theorem 2]{Williams82}.
\end{proof}

\section{An auxiliary sum and some special evaluations}

The sum $S_m(x)$ introduced in the first section can be evaluated modulo a prime through a finite hypergeometric series ${}_{2}F_1$ which depend on two parameters $\alpha$ and $\beta$ to be chosen suitably. 

\begin{theorem} Let $p$ be a prime, and let $\alpha,\beta\in\mathbb{Z}_p$ with $\beta$ prime to $p$. Then for any positive integer $m\leq \frac{p-1}{2}$,
\begin{equation}\label{Shyp}
S_m(x)\equiv \frac{1}{\beta p}\,\left(
\FF{-m}{-\beta p}{\frac{1}{2}-m-\alpha p}{\frac{1}{4x}}
-1\right)
\pmod{p}.
\end{equation}
\end{theorem}
\begin{proof}
\begin{align*}
\sum_{k=0}^{m-1}
\frac{\binom{2k}{k}x^k}{k-m}&=-\frac{1}{m}\sum_{k=0}^{m-1}\frac{\left(-m\right)_k\left(\frac{1}{2}\right)_k(4x)^k}{\left(1\right)_k\left(-m+1\right)_k}\\
&\equiv -\frac{1}{m}\sum_{k=0}^{m-1}\frac{\left(-m\right)_k\left(\frac{1}{2}+\alpha p\right)_k(4x)^k}{\left(1\right)_k\left(-m+1+\beta p\right)_k}\pmod{p}\\
&=-\frac{1}{m}\,\left(
\FF{-m}{\frac{1}{2}+\alpha p}{-m+1+\beta p}{4x}
-
\frac{\left(\frac{1}{2}+\alpha p\right)_m(-4x)^m}{\left(-m+1+\beta p\right)_m}\right)\\
&=-\frac{\left(\frac{1}{2}+\alpha p\right)_m(-4x)^m}{m\left(-m+1+\beta p\right)_m}\,\left(
\FF{-m}{-\beta p}{\frac{1}{2}-m-\alpha p}{\frac{1}{4x}}
-1\right)\\
&\equiv
\binom{2m}{m}\frac{x^m}{\beta p}\,\left(
\FF{-m}{-\beta p}{\frac{1}{2}-m-\alpha p}{\frac{1}{4x}}
-1\right)
\pmod{p}.
\end{align*}
where we used the linear transformation \cite[15.8.6]{OLBC2010}
$$\FF{-m}{b}{c}{z}=\frac{(b)_m}{(c)_m}(-z)^m
\FF{-m}{1-c-m}{1-b-m}{\frac{1}{z}}.$$
\end{proof}

Some evaluations of $S_m(x)$ are straighforward.
For $x=\frac{1}{4}$, we let $\alpha=0$ and $\beta=1$ and we find
\begin{align*}S_m\left(\frac{1}{4}\right)&\equiv \frac{1}{p}\,\left(
\FF{-m}{-p}{\frac{1}{2}-m}{1}
-1\right)=\frac{1}{p}\,\left(
\frac{\left(\frac{1}{2}-m+p\right)_m}{\left(\frac{1}{2}-m\right)_m}
-1\right)\\
&=\frac{1}{p}\,\left(
\prod_{j=1}^m\left(1-\frac{2p}{2j-1}\right)-1\right)
\equiv -2\sum_{j=1}^m\frac{1}{2j-1}=
H_m-2H_{2m}
\pmod{p},
\end{align*}
where we applied the Chu-Vandermonde identity \cite[15.4.24]{OLBC2010}
$$\FF{-m}{b}{c}{1}=\frac{(c-b)_m}{(c)_m}.$$
Another easy case is for $m=\frac{p-1}{2}$,
\begin{align*}
S_{\frac{p-1}{2}}(x)&
\equiv (-1)^n\sum_{k=0}^{n-1}\frac{\binom{n}{k}(-4)^k x^{k-n}}{k-n}\equiv -\sum_{k=0}^{n-1}\frac{\binom{n}{k} (-4x)^{k-n}}{n-k}\\
&=-\sum_{k=1}^{n}\frac{\binom{n}{k}}{k(-4x)^k}\equiv -\sum_{k=1}^{p-1}\frac{\binom{2k}{k}}{k(16x)^k}\pmod{p}.
\end{align*}
Therefore, in view of \cite{SunTauraso2010}, the following congruences hold modulo a prime $p$,
\begin{align}\label{Shalfp}
S_{\frac{p-1}{2}}\left(\frac{1}{16}\right)\equiv 0,\qquad
S_{\frac{p-1}{2}}\left(-\frac{1}{32}\right)\equiv 4q_p(2),\qquad
S_{\frac{p-1}{2}}\left(-\frac{1}{16}\right)\equiv 5f_p.
\end{align}
More evaluations can be made. 
In view of \eqref{Shyp} and by using the Maple package {\tt twoFone} written by Zeilberger \cite{Ekhad2004}, after a systematic search,
we were able to obtain various evaluations of $S_m(x)$ collected in Table 1. In order to explain how the table works we present the full proof of two rows. The other rows can be verified along the same scheme.

\begin{table}[ht]
\centering
\caption{Values of $S_m(x)$ modulo $p$ for $p\equiv 1 \pmod{r}$}\vspace{-4mm}
\renewcommand{\arraystretch}{0.5}
\begin{tabular}[t]{|c|c|c|c|c|c|}
\hline 
&&&&&\\
$r$&$m$&$x$&$\alpha,\beta$&${}_{2}F_1$&$S_m(x)$\\
&&&&&\\
\hline
&&&&&\\
3&$\frac{p-1}{6}$&$\frac{1}{16}$&
$0,\frac{1}{2}$&$(-3)^{3m}$&$q_p(3)$\\
&&&&&\\
3&$\frac{p-1}{6}$&$-\frac{1}{32}$&
$-\frac{1}{6},\frac{1}{3}$&$(-3)^{3m}$&$\frac{3}{2}q_p(3)$\\
&&&&&\\
3&$\frac{p-1}{3}$&$1$&
$-1,-1$&
$\displaystyle\frac{2^{4m}\left(\frac{3}{2}\right)_{2m}}{3^{3m}
\left(\frac{4}{3}\right)_{2m}}$&$q_p(3)$\\
&&&&&\\
3&$\frac{p-1}{3}$&$-2$&
$-1,-1$
&$\displaystyle\frac{2^{3m}\left(\frac{3}{2}\right)_{2m}}{3^{3m}
\left(\frac{5}{6}\right)_{m}\left(\frac{7}{6}\right)_{m}}$
&$-q_p(2)+\frac{3}{2}q_p(3)$\\
&&&&&\\
\hline
&&&&&\\
4&$\frac{p-1}{4}$&$-\frac{1}{4}$&$0,\frac{1}{2}$
&$\displaystyle\frac{2^{2m}\left(\frac{3}{4}\right)_{m}}{\left(\frac{1}{2}\right)_{m}}$&$\frac{3}{2}q_p(2)$
\\
&&&&&\\
\hline
&&&&&\\
5&$\frac{p-1}{10}$&$-\frac{1}{16}$&
$\frac{1}{2},\frac{1}{2}$&
$\displaystyle
\frac{5^{5m}\left(\frac{2}{5}\right)_{2m}\left(\frac{3}{5}\right)_{2m}\left(1\right)_{2m}}{\left(1\right)_{6m}}$
&$q_p(5)-2f_p$
\\
&&&&&\\
5&$\frac{p-1}{5}$&$-1$&$-1,-1$
&
$\displaystyle\frac{2^{8m}\left(\frac{3}{2}\right)_{4m}}{5^{5m}
\left(\frac{4}{5}\right)_{2m}\left(\frac{6}{5}\right)_{2m}}$
&$q_p(5)-3f_p$
\\
&&&&&\\
5&$\frac{3(p-1)}{10}$&$-\frac{1}{16}$&
$\frac{3}{2},\frac{3}{2}$
&
$\displaystyle
\frac{5^{5m}\left(\frac{4}{5}\right)_{2m}\left(\frac{6}{5}\right)_{2m}\left(1\right)_{2m}}{\left(2\right)_{6m}}$
&$q_p(5)+2f_p$\\
&&&&&\\
5&$\frac{2(p-1)}{5}$&$-1$&
$-2,-2$
&$\displaystyle\frac{2^{8m}\left(\frac{5}{2}\right)_{4m}}{5^{5m}
\left(\frac{7}{5}\right)_{2m}\left(\frac{8}{5}\right)_{2m}}$
&$q_p(5)+3f_p$\\
&&&&&\\
\hline
\end{tabular}
\end{table}

\begin{proof}[Proof of the case $p\equiv 1\pmod{5}$, $m=\frac{p-1}{5}$, $x=-1$, $\alpha=\beta=-1$]

\noindent By running {\tt twoFone} we establish the identity 
$$\FF{-m}{5m+1}{4m+\frac{3}{2}}{-\frac{1}{4}}
=\frac{2^{8m}\left(\frac{3}{2}\right)_{4m}}{5^{5m}
\left(\frac{4}{5}\right)_{2m}\left(\frac{6}{5}\right)_{2m}} 
.$$
From the basic properties of $\Gamma_p$, \eqref{de}, and Lemma \ref{LHf}, we deduce
\begin{align*}
\frac{\left(\frac{3}{2}\right)_{4m}}{
\left(\frac{4}{5}\right)_{2m}\left(\frac{6}{5}\right)_{2m}} 
&=\frac{\left(\frac{1}{2}\right)_{\frac{5m}{2}}\left(\frac{3}{2}+\frac{5m}{2}\right)_{\frac{3m}{2}}}{
\left(\frac{4}{5}\right)_{2m}\left(\frac{1}{5}\right)_{m}\left(\frac{6}{5}+m\right)_{m}}\\
&= 
\frac{\Gamma_p\left(\frac{1}{2}+\frac{5m}{2}\right)
\Gamma_p\left(\frac{3}{2}+4m\right)\Gamma_p\left(\frac{4}{5}\right)\Gamma_p\left(\frac{1}{5}\right)\Gamma_p\left(\frac{6}{5}+m\right)}{\Gamma_p\left(\frac{1}{2}\right)\Gamma_p\left(\frac{3}{2}+\frac{5m}{2}\right)\Gamma_p\left(\frac{4}{5}+2m\right)\Gamma_p\left(\frac{1}{5}+m\right)\Gamma_p\left(\frac{6}{5}+2m\right)}\\
&= 
\frac{\Gamma_p\left(\frac{7}{10}+\frac{4p}{5}\right)
\Gamma_p\left(\frac{4}{5}\right)\Gamma_p\left(\frac{1}{5}\right)}{\Gamma_p\left(\frac{1}{2}\right)\Gamma_p\left(\frac{2}{5}+\frac{2p}{5}\right)\Gamma_p\left(\frac{4}{5}+\frac{2p}{5}\right)}\\
&\equiv
\frac{\Gamma_p\left(\frac{1}{5}\right)\Gamma_p\left(\frac{7}{10}\right)
}{\Gamma_p\left(\frac{1}{2}\right)\Gamma_p\left(\frac{2}{5}\right)}
\left(1+p\left(
\frac{4}{5}H_{\frac{3(p-1)}{10}}-\frac{2}{5}H_{\frac{3(p-1)}{5}}-\frac{2}{5}H_{\frac{p-1}{5}}
\right)\right)\\
&\equiv
\frac{\Gamma_p\left(\frac{1}{5}\right)\Gamma_p\left(\frac{7}{10}\right)
}{\Gamma_p\left(\frac{1}{2}\right)\Gamma_p\left(\frac{2}{5}\right)}
\left(1+p\left(
-\frac{8q_p(2)}{5}+3f_p
\right)\right).
\end{align*}
Since by the multiplication formula \cite[p. 371]{Robert2000},
$$\Gamma_p\left(\frac{1}{5}\right)\Gamma_p\left(\frac{7}{10}\right)=\Gamma_p\left(\frac{2/5}{2}\right)
\Gamma_p\left(\frac{2/5+1}{2}\right)=
\frac{1}{2^{c}}\,\Gamma_p\left(\frac{1}{2}\right)\Gamma_p\left(\frac{2}{5}\right)$$
it follows that
$$\frac{\Gamma_p\left(\frac{1}{5}\right)\Gamma_p\left(\frac{7}{10}\right)
}{\Gamma_p\left(\frac{1}{2}\right)\Gamma_p\left(\frac{2}{5}\right)}\equiv 4^{m}\left(1-\frac{2pq_p(2)}{5}\right)
\pmod{p^2}.$$
Hence
\begin{align*}
\FF{-m}{5m+1}{4m+\frac{3}{2}}{-\frac{1}{4}}
&\equiv \frac{4^{p-1}}{5^{p-1}}
\left(1-\frac{2pq_p(2)}{5}\right)
\left(1+p\left(
-\frac{8q_p(2)}{5}+3f_p
\right)\right)\\
&\equiv 1+p(-q_p(5)+3f_p)\pmod{p^2}
\end{align*}
and by recalling that $\beta=-1$ we are done.
\end{proof}

\begin{proof}[Proof of the case $p\equiv 1\pmod{5}$, $m=\frac{3(p-1)}{10}$, $x=-\frac{1}{16}$, $\alpha=\beta=\frac{3}{2}$]

\noindent Again by  {\tt twoFone} we have the identity 
$$\FF{-m}{-5m-\frac{3}{2}}{-6m-1}{-4}
=\frac{5^{5m}\left(\frac{4}{5}\right)_{2m}\left(\frac{6}{5}\right)_{2m}\left(1\right)_{2m}}{\left(2\right)_{6m}}.$$ 
The basic properties of $\Gamma_p$, \eqref{de}, and Lemma \ref{LHf} imply
\begin{align*}
\frac{\left(\frac{4}{5}\right)_{2m}\left(\frac{6}{5}\right)_{2m}\left(1\right)_{2m}}{\left(2\right)_{6m}} 
&=\frac{\left(\frac{4}{5}\right)_{2m}\left(\frac{1}{5}\right)_{\frac{2m}{3}}\left(\frac{6}{5}+\frac{2m}{3}\right)_{\frac{4m}{3}}\left(1\right)_{2m}}{
\left(1\right)_{\frac{10m}{3}}\left(2+\frac{10m}{3}\right)_{\frac{8m}{3}}}\\
&= 
\frac{\Gamma_p\left(\frac{4}{5}+2m\right)
\Gamma_p\left(\frac{1}{5}+\frac{2m}{3}\right)\Gamma_p\left(\frac{6}{5}+2m\right)
\Gamma_p\left(1+2m\right)
\Gamma_p\left(2+\frac{10m}{3}\right)}
{\Gamma_p\left(\frac{4}{5}\right)
\Gamma_p\left(\frac{1}{5}\right)
\Gamma_p\left(\frac{6}{5}+\frac{2m}{3}\right)\Gamma_p\left(1+\frac{10m}{3}\right)
\Gamma_p\left(2+6m\right)}\\
&=
\frac{\Gamma_p\left(\frac{1}{5}+\frac{3p}{5}\right)
\Gamma_p\left(\frac{3}{5}+\frac{3p}{5}\right)\Gamma_p\left(\frac{2}{5}+\frac{3p}{5}\right)}{\Gamma_p\left(\frac{4}{5}\right)
\Gamma_p\left(\frac{1}{5}\right)
\Gamma_p\left(\frac{1}{5}+\frac{9p}{5}\right)}\\
&\equiv
\frac{\Gamma_p\left(\frac{3}{5}\right)
\Gamma_p\left(\frac{2}{5}\right)
}{\Gamma_p\left(\frac{4}{5}\right)\Gamma_p\left(\frac{1}{5}\right)}
\left(1+p\left(
\frac{3}{5}H_{\frac{2(p-1)}{5}}+\frac{3}{5}H_{\frac{3(p-1)}{5}}-\frac{6}{5}H_{\frac{4(p-1)}{5}}
\right)\right)\\
&\equiv
\frac{\Gamma_p\left(\frac{3}{5}\right)
\Gamma_p\left(\frac{2}{5}\right)
}{\Gamma_p\left(\frac{4}{5}\right)\Gamma_p\left(\frac{1}{5}\right)}
\left(1+3pf_p\right).
\end{align*}
By the reflection formula \cite[(5) p. 369]{Robert2000},
$$\frac{\Gamma_p\left(\frac{3}{5}\right)
\Gamma_p\left(\frac{2}{5}\right)
}{\Gamma_p\left(\frac{4}{5}\right)\Gamma_p\left(\frac{1}{5}\right)}=1,$$
and therefore
\begin{align*}
\FF{-m}{-5m-\frac{3}{2}}{-6m-1}{-4}
&\equiv 5^{\frac{3(p-1)}{2}}
\left(1+3pf_p\right)\equiv 
\left(1+p\frac{3q_p(5)}{2}\right)\left(1+3pf_p\right)\\
&\equiv \left(1+p\left(\frac{3q_p(5)}{2}+3f_p\right)\right)\pmod{p^2}.
\end{align*}
The desired conclusion follows as soon as we note that $\beta=\frac{3}{2}$.
\end{proof}

\section{Proof of the main theorem}

This opening lemma introduces two combinatorial identities.

\begin{lemma} For any $1\leq m\leq n$, and for any $z\in\mathbb{C}$,
\begin{equation}\label{polid1}
\left(1+\frac{1}{z}\right)^{n+1}\sum_{k=1}^{m}\frac{z^{k}}{k\binom{n}{k}}
=
H_n+H_m-H_{m'}-\left(T_n(z)+T_m(1/z)-T_{m'}(z)\right)
\end{equation}
and
\begin{equation}\label{polid2}
\left(1+\frac{1}{z}\right)^{n+1}\sum_{k=m+1}^{n}\frac{z^{k}}{k\binom{n}{k}}
=
H_n+H_{m'}-H_m-\left(T_n(1/z)+T_{m'}(z)-T_m(1/z)\right)
\end{equation}
where $m'=n-m$, and
$$T_r(z)=\frac{1}{\binom{n}{r}}\sum_{k=0}^{r-1}\binom{n}{k}\frac{z^{k-r}}{k-r}.
$$
\end{lemma}
\begin{proof} After multiplying \eqref{polid1} by $z^{n}$, we find
$$
(z+1)^{n+1}\sum_{k=1}^{m}\frac{z^{k-1}}{k\binom{n}{k}}
=
\sum_{k=0}^{n-1}\binom{n}{k}\frac{z^{k}}{n-k}
+(H_n+H_m-H_{n-m})z^n
-\frac{1}{\binom{n}{m}}\sum_{\substack{k=0\\k\neq m}}^{n}\binom{n}{k}\frac{z^{n+m-k}}{k-m}.
$$
Both sides of  are polynomials of degree $n+m$, so it suffices to compare the coefficients of $z^d$ for each $d$ such that $0\leq d\leq n+m$:
$$\sum_{k=1}^{m}\frac{\binom{n+1}{d-k+1}}{k\binom{n}{k}}=
\begin{cases}\displaystyle
\frac{1}{n-d}\left(
\binom{n}{d}-\frac{\binom{n}{n+m-d}}{\binom{n}{m}}\right)
&\text{if $d\not=n$},\\\\
H_n+H_m-H_{n-m} &\text{if $d=n$,}
\end{cases}
$$
which can be easily verified by induction with respect to $m$.
As regards \eqref{polid2}, just subtract \eqref{polid1} from \eqref{polid1} where we set $m=n$.
\end{proof}

Now we are able to show that the sum $S_m(x)$ is truly relevant to our needs.

\begin{theorem}
Let $p$ be an odd prime with $1\leq m<\frac{p-1}{2}$. Let $x\in\mathbb{Z}_p$ be such that $x(1-4x)$ is prime to $p$. Let
\begin{equation}\label{defQ}
Q_p(m,x)=H_{\frac{p-1}{2}}+H_m-H_{m'}-\left(S_{\frac{p-1}{2}}(y)+S_m(x)-S_{m'}(y)\right)
\end{equation}
where $m'=\frac{p-1}{2}-m$ and $y=\frac{1}{16x}$.
Then
\begin{align}\label{fcong1}
\sum_{k=p-m}^{p-1}\binom{2k}{k}x^k
&\equiv
-\frac{2px Q_p(m,x)}{(1-4x)^{\frac{p+1}{2}}}\pmod{p^2}
\end{align}
and
\begin{align}\label{fcong2}
\sum_{k=\frac{p+1}{2}}^{p-1-m'}\binom{2k}{k}y^k&\equiv
-\frac{2py Q_p(m,x)}{(1-4y)^{\frac{p+1}{2}}}\pmod{p^2}.
\end{align}
\end{theorem}
\begin{proof}
We first note that for $1\leq k\leq n=\frac{p-1}{2}$, 
$$\binom{n}{k}\equiv \binom{2k}{k}(-4)^{-k}\pmod{p}.$$
Hence for $z=\frac{1}{4x}$, we find
$$T_r(z)\equiv S_r\left(-\frac{1}{4x}\right)  \pmod{p}$$
and, by \eqref{polid1}, 
\begin{align*}
\left(1-4x\right)^{\frac{p+1}{2}}\sum_{k=1}^{m}\frac{(1/x)^{k}}{k\binom{2k}{k}}
\equiv
Q_p(x,m)\pmod{p}.
\end{align*}
Futhermore, for $1\leq k\leq \frac{p-1}{2}$, 
$$k\binom{2k}k\binom{2(p-k)}{p-k}\equiv -2p\pmod{p^2}$$
which imply
\begin{align*}
\sum_{k=p-m}^{p-1}\binom{2k}{k}x^k
&=\sum_{k=1}^{m}\binom{2(p-k)}{p-k}x^{p-k}
\equiv -2px^p\sum_{k=1}^{m}\frac{(1/x)^k}{k\binom{2k}{k}}\pmod{p^2},
\end{align*}
and \eqref{fcong1} follows immediately. In a similar way, by applying  \eqref{polid2}, we show \eqref{fcong2}.
\end{proof}

Finally, all is set up for our final proof.

\begin{proof}[Proof of Theorem \ref{Thadam}] If $p\equiv 1$ modulo $3$, by \eqref{defQ}, Lemma \ref{LHf}, \eqref{Shalfp} and Table 1,  
\begin{align*}
Q_p\left(\frac{p-1}{3},-2\right)&=
H_{\frac{p-1}{2}}+H_{\lfloor\frac{p}3\rfloor}-H_{\lfloor\frac{p}6\rfloor}-S_{\frac{p-1}{2}}\left(-\frac{1}{32}\right)-S_{\frac{p-1}{3}}(-2)+S_{\frac{p-1}{6}}\left(-\frac{1}{32}\right)\\
&\equiv
-2q_p(2)-\frac{3q_p(3)}{2}+2q_p(2)+\frac{3q_p(3)}{2}
-4q_p(2)+q_2(p)-\frac{3q_p(3)}{2}+\frac{3q_p(3)}{2}
\\
&\equiv -3q_p(2) \pmod{p}.
\end{align*}
Thus, by \eqref{fcong1} and \eqref{fcong2} together with \eqref{ncong1a} and \eqref{ncong1b}, we obtain the first pair \eqref{pair1},
\begin{align*}
&\sum_{k=1}^{\lfloor\frac{2p}3\rfloor}\binom{2k}{k}(-2)^k
\equiv \sum_{k=1}^{p-1}\binom{2k}{k}(-2)^k-\frac{4p Q_p\left(\frac{p-1}{3},-2\right)}{9^{\frac{p+1}{2}}}\equiv 0 \pmod{p^2},\\
&\sum_{k=0}^{\lfloor\frac{5p}6\rfloor}\frac{\binom{2k}{k}}{(-32)^k}
\equiv \sum_{k=0}^{\frac{p-1}{2}}\frac{\binom{2k}{k}}{(-32)^k}
+\frac{p Q_p\left(\frac{p-1}{3},-2\right)}{16(1+1/8)^{\frac{p+1}{2}}}
\equiv \left(\frac{2}{p}\right) \pmod{p^2}.
\end{align*}
Moreover, if $p\equiv 1$ modulo $5$ then by \eqref{defQ}, Lemma \ref{LHf}, \eqref{Shalfp} and Table 1,  
\begin{align*}
Q_p\left(\frac{p-1}{5},-1\right)&=
H_{\frac{p-1}{2}}+H_{\lfloor\frac{p}5\rfloor}-H_{\lfloor\frac{3p}{10}\rfloor}
-S_{\frac{p-1}{2}}\left(-\frac{1}{16}\right)
-S_{\frac{p-1}{5}}(-1)
+S_{\frac{3(p-1)}{10}}\left(-\frac{1}{16}\right)\\
&\equiv
-2q_p(2)-\frac{5q_p(5)}{4}-\frac{5f_p}{4}+2q_p(2)+\frac{5q_p(5)}{4}-\frac{15f_p}{4}\\
&\qquad -5f_p-q_p(5)+3f_p+q_p(5)+2f_p\equiv -5f_p \pmod{p},
\end{align*}
and 
\begin{align*}
Q_p\left(\frac{2(p-1)}{5},-1\right)&=
H_{\frac{p-1}{2}}+H_{\lfloor\frac{2p}5\rfloor}-H_{\lfloor\frac{p}{10}\rfloor}
-S_{\frac{p-1}{2}}\left(-\frac{1}{16}\right)
-S_{\frac{2(p-1)}{5}}(-1)
+S_{\frac{p-1}{10}}\left(-\frac{1}{16}\right)\\
&\equiv
-2q_p(2)-\frac{5q_p(5)}{4}+\frac{5f_p}{4}+2q_p(2)+\frac{5q_p(5)}{4}+\frac{15f_p}{4}
\\
&\qquad -5f_p-q_p(5)-3f_p+q_p(5)-2f_p\equiv -5f_p \pmod{p}.
\end{align*}
Then the second pair \eqref{pair2} and the third pair \eqref{pair3} come from \eqref{fcong1}, \eqref{fcong2}, \eqref{ncong2a} 
 and \eqref{ncong2b}.
\end{proof}

\vskip 3mm \noindent{\bf Acknowledgments.}
The first author is funded by the Startup Foundation for Introducing Talent of Nanjing University of Information Science and Technology (2019r062).

\end{document}